\documentclass[11pt]{article}%
\usepackage{amsmath}
\usepackage{amsfonts}
\usepackage{amssymb}
\usepackage{graphicx}
\usepackage{subfigure}%
\newtheorem{theorem}{Theorem}

\newtheorem{corollary}[theorem]{Corollary}

\newtheorem{definition}[theorem]{Definition}

\newtheorem{lemma}[theorem]{Lemma}

\newenvironment{proof}[1][Proof]{\noindent\textbf{#1.} }{\ \rule{0.5em}{0.5em}}
\begin{document}

\title{How to Regularize a Symplectic-Energy-Momentum Integrator}
\author{Yosi Shibberu\\Mathematics Department\\Rose-Hulman Institute of Technology\\Terre Haute, IN 47803\\shibberu@rose-hulman.edu\\http://rose-hulman.edu/$\sim$shibberu/DTH\_Dynamics/DTH\_Dynamics.htm}
\maketitle

\begin{abstract}
We identify ghost trajectories of symplectic-energy-momentum (SEM)
integration and show that the ghost trajectories are not time
reversible. We explain how SEM integration can be regularized, in
a SEM preserving manner, so that it is time reversible. We
describe an algorithm for implementing the regularized SEM
integrator. Simulation results for the pendulum are given.
Coordinate invariance of the regularized SEM integrator is briefly
discussed.

\begin{description}
\item[Key Words] DTH dynamics, symplectic energy momentum integrator, discrete
mechanics, discrete time Hamiltonian, discrete variational
principles, principle of least action, energy conserving methods,
extended phase space, midpoint method.

\end{description}
\end{abstract}

\section{Introduction}

Is symplectic-energy-momentum (SEM)\ integration obstructed by singularities?
In 1994, the answer appeared to be yes. SEM integration of the rotations of a
pendulum did not appear possible, Shibberu \cite{Shibberu-94}. In this
article, we show that SEM integration is actually \emph{not} obstructed by
singularities. Difficult to compute, ghost trajectories are identified, but
these ghost trajectories are not time reversible. We explain how SEM
integration can be regularized, \emph{in a SEM conserving manner, }so that it
is time reversible. We begin with a brief review of SEM integration.

A SEM integrator is a symplectic integrator which exactly conserves energy and
momentum. Symplecticness implies the integrator can be derived from a discrete
variational principle \cite{Wu-90}. The discrete variational principle of a
symplectic integrator gives it coordinate invariant properties.

SEM integration emerged from two lines of research, symplectic
integration and discrete mechanics. Efficient computation is
emphasized in symplectic integration. Preservation of the physical
laws of nature is the emphasis in discrete mechanics. The term
\textquotedblleft symplectic-energy-momentum
integrator\textquotedblright\ was coined and popularized by Kane,
Marsden and Ortiz \cite{Kane-99}. See also Chen, Guo and Wu
\cite{Chen-03} for related work on higher-order, sympletic-energy
integrators. Guibout and Bloch \cite{Guibout-04} have developed a
general framework for deriving many of the published symplectic
integrators, including SEM integrators. They also provide an
interesting comparison of the various discrete variational
principles used.

The author's work on a discrete-time theory for Hamiltonian
dynamics (DTH dynamics) \cite{Shibberu-92} predates the work of
Kane et. al. \cite{Kane-99}. DTH dynamics is a SEM integrator. DTH
dynamics originated from an effort to obtain the exact energy and
momentum conserving properties of the discrete mechanics of
Greenspan \cite{Greenspan-74}, \cite{Greenspan-80} from the
variational principle used in the discrete mechanics of Lee
\cite{Lee-83}, \cite{Lee-87}. DTH dynamics was proved in 1994 (see
Shibberu \cite{Shibberu-94}, \cite{Shibberu-97}) to be symplectic
and hence a SEM integrator. D'Innocenzo, Renna and Rotelli
\cite{DInnocenzo-87} have done work that can also be related to
SEM integration.

In the author's work, SEM integration is accomplished by varing the time step
of the midpoint scheme to enforce exact energy conservation at the
\emph{midpoint} of each step. Symplecticness and momentum conservation occur
at the \emph{vertices} of each step. The relationship between the time step
and energy conservation originates from the fact that the negative of the
energy (Hamiltonian) is the momentum corresponding to time, Shibberu
\cite{Shibberu-92}, Lee \cite{Lee-83}.

The requirements of symplectic-energy integration are highly
restrictive as illustrated by Ge's Theorem \cite{Ge-88},
\cite{Ge-91}. An early existence and uniqueness result for DTH
dynamics was given in Shibberu \cite{Shibberu-92} and an
explanation of why Ge's Theorem is not violated was given in
Shibberu \cite{Shibberu-97}. The sufficient condition for the
existence and uniqueness of DTH trajectories proved in Shibberu
\cite{Shibberu-92} does not cover all the points in phase space
where the Hamiltonian function is smoothly defined. It first
appeared, from simulation results in Shibberu \cite{Shibberu-94},
that DTH dynamics was obstructed by points where this sufficient
condition does not hold. Kane et. al. \cite{Kane-99} later
observed similar difficulties near
points they refer to as \textquotedblleft turning points\textquotedblright%
\ and Chen et. al. \cite{Chen-03} also mentioned the need to avoid
singularities in their algorithm. In this article, we explain how
SEM integration can be regularized in a manner which preserves SEM
properties and time reversibility at such points.

The outline of this paper is as follows. In section
\ref{foundations}, we review the foundations of DTH dynamics. We
introduce a discretization of Hamiltonian dynamics which is
equivalent to, but simpler than the discretization used in
Shibberu \cite{Shibberu-97}. An example of a ghost trajectory of
the pendulum is also given. In section \ref{regularization}, we
introduce the two complementary variational principles used to
regularized DTH dynamics. Regularized DTH equations are derived
and shown to preserve SEM properties. Coordinate invariance of the
regularized equations is briefly discussed. In section
\ref{numerics}, we give a detailed description of an algorithm for
solving the regularized DTH equations. Numerical results for the
pendulum and Kepler's one-body problem are discussed. Finally,
certain peculiarities of regularized DTH dynamics are described.

\section{DTH Dynamics\label{foundations}}

\subsection{Foundations of DTH Dynamics}

We begin by introducing an extended-phase space version of the principle of
least action. Let $H(t,q_{1},\ldots,q_{n},p_{1},\ldots,p_{n})$ be the
Hamiltonian function of an $n$-degree of freedom Hamiltonian dynamical system
where $t$ is time and $q_{i}$ and $p_{i},$ $i=1,\ldots,n,$ are position and
momentum coordinates. Let $q=\left(  q_{1},\ldots,q_{n},t\right)  ^{\top}$ and
$p=\left(  p_{1},\ldots,p_{n},\wp\right)  ^{\top}$ be extended phase space
coordinates where $\wp$ is the momentum conjugate to time. (See
\cite{Lanczos-70}, \cite{Goldstein-80} or \cite{Shibberu-94} for a description
of $\wp.$) Let $z=(q,p)^{\top}.$ Consider the extended-phase space action
integral%
\[
\mathcal{A}(z(\cdot))=\int_{\tau_{0}}^{\tau_{f}}\frac{1}{2}z(\tau)^{\top
}Jz^{\prime}(\tau)^{\top}\,d\tau,\ \text{where }J=\left(
\begin{array}
[c]{rr}%
0 & I\\
-I & 0
\end{array}
\right)
\]
and $I$ is the $n+1$ dimensional identity matrix. The extended-phase space
Hamiltonian function is $\mathcal{H}(z)=\wp+H(t,q_{1},\ldots,q_{n}%
,p_{1},\ldots,p_{n}).$ The principle of least (stationary) action states that
the trajectory $z(\tau)$ of a Hamiltonian dynamical system cause the action
integral $\mathcal{A}(z(\cdot))$ to be stationary under variations which
satisfy the boundary constraints $q(\tau_{0})=q_{0},$ $p(\tau_{f})=p_{f}$ and
the Hamiltonian constraint $\mathcal{H}(z)\equiv0.$ (Additional details are
given in \cite{Shibberu-97}.)

The action integral $\mathcal{A}(z(\cdot))$ is discretized in
\cite{Shibberu-97} by evaluating the integral along piecewise-linear,
continuous trajectories and then appending boundary terms to account for the
boundary conditions. An equivalent discretization with no boundary terms is
described below. This discretization makes it possible to provide a simpler
derivation of the DTH equations.%

%TCIMACRO{\TeXButton{Figure: triangle}{\begin{figure}
%\begin{center}
%\begin{picture}(100,100)
%\put(10,20){\circle{5}} \put(10,20){\line(1,0){25}}
%\put(35,20){\circle{5}} \put(32,9){$z_k$}
%\put(35,20){\line(1,1){30}} \put(47,35){\line(1,0){5}}
%\put(58,29){$\overline{z}_{k}$} \put(65,50){\circle{5}}
%\put(65,50){\line(1,0){25}} \put(60,57){$z_{k+1}$}
%\put(90,50){\circle{5}} \multiput(35,20)(0,5){6}{\line(0,1){2}}
%\multiput(65,50)(-5,0){7}{\line(1,0){2}} \put(39,41){$\sigma_k$}
%\put(-5,56){$(q_k,p_{k+1})$}
%\end{picture}
%\end{center}
%\caption{A piecewise-linear, continuous trajectory in
%extended-phase space.} \label{triangle}
%\end{figure} }}%
%BeginExpansion
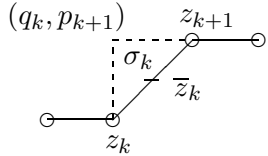
\begin{figure}
\begin{center}
\begin{picture}(100,100)
\put(10,20){\circle{5}} \put(10,20){\line(1,0){25}}
\put(35,20){\circle{5}} \put(32,9){$z_k$}
\put(35,20){\line(1,1){30}} \put(47,35){\line(1,0){5}}
\put(58,29){$\overline{z}_{k}$} \put(65,50){\circle{5}}
\put(65,50){\line(1,0){25}} \put(60,57){$z_{k+1}$}
\put(90,50){\circle{5}} \multiput(35,20)(0,5){6}{\line(0,1){2}}
\multiput(65,50)(-5,0){7}{\line(1,0){2}} \put(39,41){$\sigma_k$}
\put(-5,56){$(q_k,p_{k+1})$}
\end{picture}
\end{center}
\caption{A piecewise-linear, continuous trajectory in
extended-phase space.} \label{triangle}
\end{figure}
%EndExpansion

\begin{lemma}
\label{action_triangle}Let $\partial\sigma_{k}$ represent the boundary of the
triangle $\sigma_{k}$ in extended-phase space shown in Figure \ref{triangle}.
Then, along $\partial\sigma_{k},$
\[
\mathcal{A}(z(\cdot))=\int_{\partial\sigma_{k}}\frac{1}{2}z(\tau)^{\top
}Jz^{\prime}(\tau)^{\top}\,d\tau=\frac{1}{2}\Delta q_{k}{}^{\top}\Delta
p_{k}.
\]

\end{lemma}

Lemma \ref{action_triangle} follows from Stoke's formula \cite{Arnold-78}.

\begin{definition}
[One Step Action]\label{one_step_action}The one-step action of a discrete-time
Hamiltonian dynamical system is defined to be%
\[
\mathcal{A}(z_{k},z_{k+1})=\frac{1}{2}\Delta q_{k}{}^{\top}\Delta
p_{k},\ \text{where }z_{k}=(q_{k},p_{k})^{\top}.
\]

\end{definition}

A discrete-time Hamiltonian (DTH) trajectory is a piecewise-linear, continuous
trajectory which satisfies the following discrete variational principle.

\begin{definition}
[DTH Principle of Stationary Action]\label{DTH_principle}The one-step action
$\mathcal{A}(z_{k},z_{k+1})$ is stationary along a DTH trajectory for
variations which fix $q_{k}\ $and $p_{k+1}$ and satisfy the Hamiltonian
constraint $\mathcal{H}(\overline{z}_{k})=0$ where $\overline{z}_{k}=\frac
{1}{2}(z_{k+1}+z_{k})$ and $k=0,\ldots,N-1.$
\end{definition}

\begin{theorem}
[DTH Equations]\label{DTH_equations}A DTH trajectory is determined by the
following equations:%
\begin{subequations}
\begin{align}
\Delta z_{k}  &  =\lambda_{k}J\mathcal{H}_{z}(\overline{z}_{k}%
)\label{DTH_equations_1}\\
\mathcal{H}(\overline{z}_{k})  &  =0. \label{DTH_equations_2}%
\end{align}

\end{subequations}
\end{theorem}

The proof of Theorem \ref{DTH_equations} follows from the proof of Theorem
\ref{reg_DTH_equations}. The proof that equations (\ref{DTH_equations_1}%
)--(\ref{DTH_equations_2}) is a SEM integrator (i.e. preserve symplecticness
and conserve momentum at $z_{k}$ and conserves energy at $\overline{z}_{k}$ )
can be found in \cite{Shibberu-97}.

\subsection{An Example of SEM Integration}

The extended-phase space Hamiltonian function for a pendulum is
$\mathcal{H}(q,p,\wp)=\wp+\frac{1}{2}p^{2}-\cos(q).$ The
corresponding DTH
equations are%
\begin{subequations}
\begin{align}
\Delta q_{k}  &  =\lambda_{k}\overline{p}_{k}\label{DTH_pendulum_1}\\
\Delta t_{k}  &  =\lambda_{k}\label{DTH_pendulum_2}\\
\Delta p_{k}  &  =-\lambda_{k}\sin(\overline{q}_{k})\label{DTH_pendulum_3}\\
\Delta\wp_{k}  &  =0\label{DTH_pendulum_4}\\
0  &  =\overline{\wp}_{k}+\frac{1}{2}\overline{p}_{k}^{2}-\cos(\overline
{q}_{k}) \label{DTH_pendulum_5}%
\end{align}
Observe from equation (\ref{DTH_pendulum_2}) that $\lambda_{k}$
equals the time step $\Delta t_{k}$.

Figure \ref{phase_portrait} shows two DTH trajectories projected
onto the phase portrait of the pendulum. Observe that the linear
segments of DTH trajectories are tangent to their respective
energy-conserving manifolds (except possibly where they cross the v-shaped curves).%

A sufficient condition for the existence and (local) uniqueness of solutions
to equations (\ref{DTH_pendulum_1})--(\ref{DTH_pendulum_5}) is the condition
$\psi(z)\neq0$ where $\psi(z)=(J\mathcal{H}_{z})^{\top}\mathcal{H}%
_{zz}(J\mathcal{H}_{z})$ \cite{Shibberu-92}. The function
$\psi(z)$ plays a key role in the regularization of SEM
integration. The v-shaped curves in Figure \ref{phase_portrait}
are points in phase space where $\psi(z)=0$ and were first
described in \cite{Shibberu-94}.

%TCIMACRO{\TeXButton{Figure: phase_portrait}{\begin{figure}
%\begin{center}
%\subfigure[Phase
%portrait]{\label{phase_portrait}\includegraphics
%[height=1.5in]{phase_portrait.eps}}
%\subfigure[Condition
%number]{\label{condition_number}\includegraphics
%[height=1.5in]{condition_number.eps}}
%\end{center}
%\caption{DTH trajectories for the nonlinear pendulum. The v-shaped
%curves correspond to points where $\psi=0$.} \label{DTH_trj}
%\end{figure}}}%
%BeginExpansion
\begin{figure}
\begin{center}
\subfigure[Phase
portrait]{\label{phase_portrait}\includegraphics
[height=1.5in]{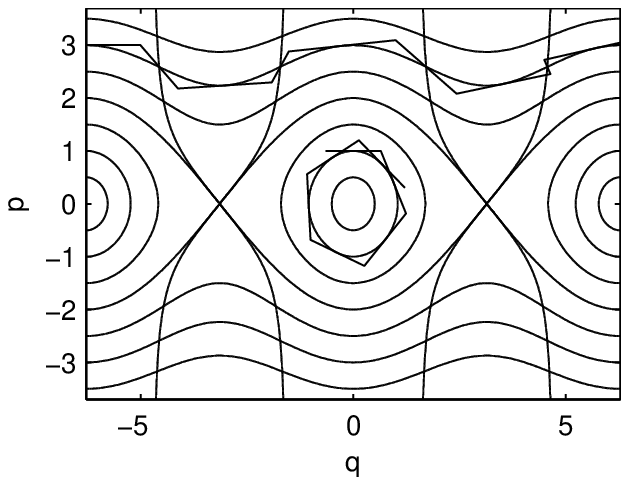}}
\subfigure[Condition
number]{\label{condition_number}\includegraphics
[height=1.5in]{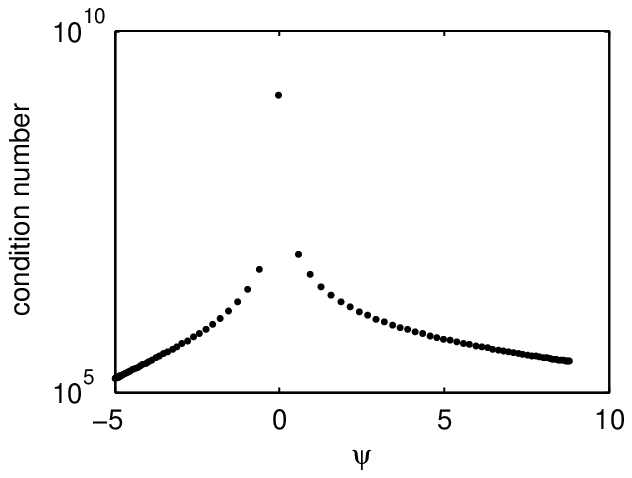}}
\end{center}
\caption{DTH trajectories for the nonlinear pendulum. The v-shaped
curves correspond to points where $\psi=0$.} \label{DTH_trj}
\end{figure}%
%EndExpansion

\subsection{Ghost Trajectories}

Figure \ref{condition_number} shows that the DTH equations
(\ref{DTH_pendulum_1})--(\ref{DTH_pendulum_5}) are poorly conditioned when
$\psi\ $is near zero. Never the less, these equations determine trajectories
which cross the v-shaped curves in Figure \ref{phase_portrait} in a SEM
conserving manner. However, these trajectories are not time reversible as is
seen in Figure \ref{irreversible}. The linear segment crossing $\psi(z)=0$
forward in time is tangent to the energy-conserving manifold on the left side,
but the segment crossing backward in time is tangent on the opposite side. We
will call these non-reversible trajectories, ghost trajectories. The
trajectory shown crossing the v-shaped curves in Figure \ref{phase_portrait}
has been regularized in a SEM conserving manner so that it is time reversible.
How this is done is explained in the next section.%

%TCIMACRO{\TeXButton{Figure: irreversible}{\begin{figure}
%\begin{center}
%\subfigure[View of $\cal{H}$ = 0, $\psi$ =
%0]{\label{irreversibility}\includegraphics[height=1.5in]{irreversibility.eps}}
%\subfigure[Magnified view of
%$\cal{H}$ = 0, $\psi$ = 0]{\label{irreversibility_2}\includegraphics
%[height=1.5in]{irreversibility_2.eps}}
%\end{center}
%\caption{A ghost DTH trajectory crossing $\psi= 0$ forward in
%time (solid curve) and then time-reversed so that it crosses
%backward in time (dashed curve).} \label{irreversible}
%\end{figure}}}%
%BeginExpansion
\begin{figure}
\begin{center}
\subfigure[View of $\cal{H}$ = 0, $\psi$ =
0]{\label{irreversibility}\includegraphics[height=1.5in]{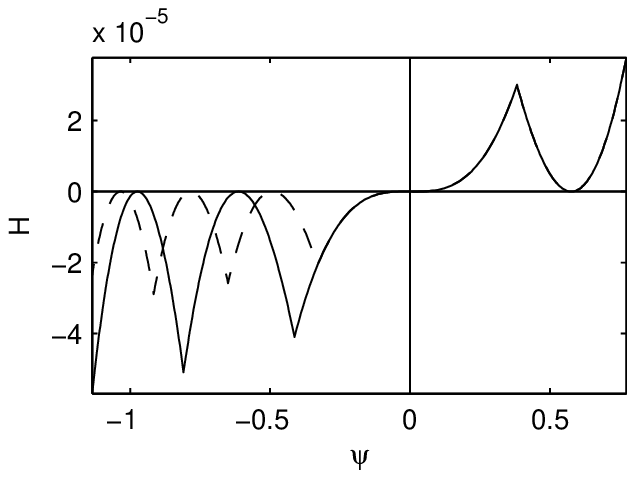}}
\subfigure[Magnified view of
$\cal{H}$ = 0, $\psi$ = 0]{\label{irreversibility_2}\includegraphics
[height=1.5in]{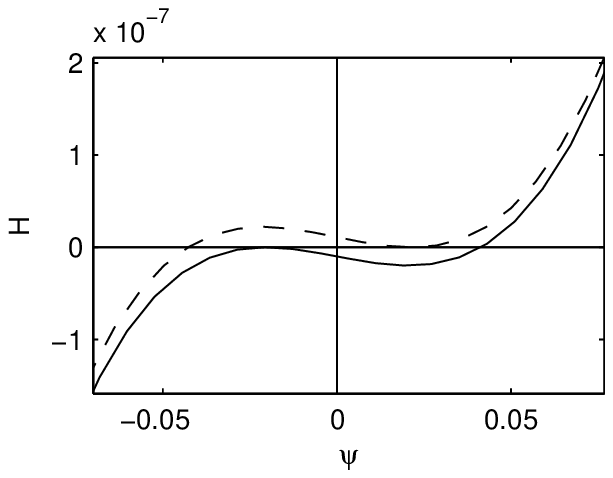}}
\end{center}
\caption{A ghost DTH trajectory crossing $\psi= 0$ forward in
time (solid curve) and then time-reversed so that it crosses
backward in time (dashed curve).} \label{irreversible}
\end{figure}%
%EndExpansion

\section{Regularized DTH Dynamics\label{regularization}}

\subsection{Regularization}

We regularize DTH dynamics by resorting to two complementary variational
principles. The first variational principle restricts variations in Definition
\ref{DTH_principle} by the inequality constraint $\psi(\overline{z}_{k}%
)\geq\psi_{k}$ where $\psi_{k}$ is a constant. The second variational
principle restricted variations by the inequality constraint $\psi
(\overline{z}_{k})\leq\psi_{k}$. We alternate between the two variational
principles to generate trajectories which cross $\psi=0$ in a time reversible manner.
\end{subequations}
\begin{definition}
[Regularized DTH Principle]The one-step action $\mathcal{A}(z_{k},z_{k+1})$ is
stationary along a DTH trajectory for variations which fix $q_{k}\ $and
$p_{k+1}$ and satisfy the Hamiltonian constraint $\mathcal{H}(\overline{z}%
_{k})=0$ and the inequality constraint $\psi(\overline{z}_{k})\geq\psi_{k}$
(or $\psi(\overline{z}_{k})\leq\psi_{k}$).
\end{definition}

\begin{theorem}
[Regularized DTH Equations]\label{reg_DTH_equations} Assume $\mathcal{H}%
_{z}(\overline{z}_{k})$ and $\psi_{z}(\overline{z}_{k})$ are linearly
independent for $k=0,\ldots,N-1.$ A regularized DTH trajectory must satisfy
the following equations and inequalities:%
\begin{subequations}
\begin{align}
\Delta z_{k}  &  =\lambda_{k}J\mathcal{H}_{z}(\overline{z}_{k})+\mu_{k}%
J\psi_{z}(\overline{z}_{k})\label{reg_DTH_eq_1}\\
\mathcal{H}(\overline{z}_{k})  &  =0\label{reg_DTH_eq_2}\\
\psi(\overline{z}_{k})  &  \geq\psi_{k}\text{\ (or\ }\psi(\overline{z}%
_{k})\leq\psi_{k}\text{)}\label{reg_DTH_eq_3}\\
\mu_{k}\left(  \psi(\overline{z}_{k})-\psi_{k}\right)   &
=0\label{reg_DTH_eq_4}\\
\mu_{k}  &  \leq0\text{\ (or\ }\mu_{k}\geq0\text{)} \label{reg_DTH_eq_5}%
\end{align}

\end{subequations}
\end{theorem}

\begin{proof}
Define the Lagrangian function
\[
\mathcal{L(}z_{k},z_{k+1},\lambda_{k},\mu_{k})=\mathcal{A(}z_{k}%
,z_{k+1})+\lambda_{k}\mathcal{H(}\overline{z}_{k})+\mu_{k}\mathcal{\psi
(}\overline{z}_{k})
\]
where $\lambda_{k}$ is a Lagrange multiplier for the equality constraint
$\mathcal{H(}\overline{z}_{k})=0$ and $\mu_{k}$ is a Karush-Kuhn-Tucker (KKT)
multiplier for the inequality constraint $\psi(\overline{z}_{k})\geq\psi_{k}%
$\ (or\ $\psi(\overline{z}_{k})\leq\psi_{k}$). Applying the KKT necessary
conditions, \cite{Bazaraa-93}, \cite{Chong-01}, to the regularized DTH
principle results in the following equations:%
\begin{subequations}
\begin{align}
\mathcal{L}_{p_{k}}  &  =-\frac{1}{2}\Delta q_{k}+\frac{1}{2}\lambda
_{k}\mathcal{H}_{p}\mathcal{(}\overline{z}_{k})+\frac{1}{2}\mu_{k}%
\mathcal{\psi}_{p}\mathcal{(}\overline{z}_{k})=0\label{Lp_eq}\\
\mathcal{L}_{q_{k+1}}  &  =\frac{1}{2}\Delta p_{k}+\frac{1}{2}%
\lambda_{k}\mathcal{H}_{q}\mathcal{(}\overline{z}_{k})+\frac{1}{2}\mu
_{k}\mathcal{\psi}_{q}\mathcal{(}\overline{z}_{k})=0\label{Lq_eq}\\
\mathcal{L}_{\lambda_{k}}  &  =\mathcal{H}(\overline{z}_{k})=0\nonumber\\
\psi(\overline{z}_{k})  &  \geq\psi_{k}\ (\text{or\ }\psi(\overline{z}%
_{k})\leq\psi_{k})\nonumber\\
\mu_{k}\left(  \psi(\overline{z}_{k})-\psi_{k}\right)   &  =0\nonumber\\
\mu_{k}  &  \leq0\text{,\ (or\ }\mu_{k}\geq0\text{).}\nonumber
\end{align}
Equations (\ref{Lp_eq}) and (\ref{Lq_eq}) can be rearranged and combined to
give equation (\ref{reg_DTH_eq_1}).
\end{subequations}
\end{proof}

Time reversibility of the regularized DTH trajectory follows from the
following observation. If the inequality constraint $\psi(\overline{z}%
_{k})\geq\psi_{k}$ is active forward in time, the inequality constraint
$\psi(\overline{z}_{k})\leq\psi_{k}$ is active backward in time. Therefore,
the same equation, $\psi(\overline{z}_{k})=\psi_{k},$ applies for both
directions in time.

\subsection{Symplectic-Energy-Momentum Properties}

Next, we show that the regularized DTH equations (\ref{reg_DTH_eq_1}%
)--(\ref{reg_DTH_eq_5}) is a SEM integrator. Conservation of energy
(Hamiltonian) follows directly from equation (\ref{reg_DTH_eq_2}). To prove
symplecticness, we identify a generating function which generates symplectic
transformations between adjacent vertices of a regularized DTH trajectory.

\begin{theorem}
[Generating Function]\label{generating_function}Assume $\mu_{k}$ and
$\mathcal{\psi(}\overline{z}_{k})$ are not simultaneous equal to zero. Then
the function%
\[
S(q_{k},p_{k+1})=q_{k}{}^{\top}p_{k+1}+\mathcal{L(}z_{k},z_{k+1},\lambda
_{k},\mu_{k})\
\]
is a generating function which determines a symplectic transformation between
the vertices $z_{k}$ and $z_{k+1}$ of a regularized DTH trajectory for
$k=0,\ldots,N-1.$ The function $\mathcal{L(}z_{k},z_{k+1},\lambda_{k},\mu
_{k})$ is the Lagrangian function used in the proof of Theorem
\ref{reg_DTH_equations} . The variables $q_{k+1},\ p_{k},\ \lambda_{k}$ and
$\mu_{k}$ are determined by the regularized DTH equations (\ref{reg_DTH_eq_1}%
)--(\ref{reg_DTH_eq_5}).
\end{theorem}

\begin{proof}
We show that $S_{q_{k}}(q_{k},p_{k+1})=p_{k}.$ The equation $S_{p_{k+1}}%
(q_{k},p_{k+1})=q_{k}$ follows in a similar fashion.
\begin{align*}
S_{q_{k}}  &  =p_{k+1}+\mathcal{L}_{q_{k}}\\
&  =p_{k+1}\mathcal{-}\frac{1}{2}\Delta p_{k}+\frac{1}{2}\lambda
_{k}\mathcal{H}_{q}\mathcal{(}\overline{z}_{k})+\frac{1}{2}\mu_{k}%
\mathcal{\psi}_{q}\mathcal{(}\overline{z}_{k})\\
&  =p_{k+1}\mathcal{-}\frac{1}{2}\Delta p_{k}\mathcal{-}\frac{1}{2}\Delta
p_{k}\\
&  =p_{k}%
\end{align*}

\end{proof}

We note that the transformation $z_{k+1}(z_{k})$ may not be differentiable
when both $\mathcal{\psi(}\overline{z}_{k})\ $and $\mu_{k}\ $are equal to
zero. Theorem \ref{generating_function} may not be valid for this special case.

The following lemma will be used to prove conservation of momentum.
(Conservation of momentum is restricted here to linear and quadratic functions
e.g. linear and angular momentum in Cartesian coordinates.)

\begin{lemma}
\label{psi_bracket}If $L(z)$ is a quadratic function and the Poisson bracket
$[L,\mathcal{H}]$ is identically equal to zero, then the Poisson bracket
$[L,\psi]$ is identically equal to zero.
\end{lemma}

\begin{proof}
By assumption, $[L,\mathcal{H}]=L_{z}^{\top}J\mathcal{H}_{z}=J^{i_{1}i_{2}%
}L_{i_{1}}\mathcal{H}_{i_{2}}=0.$ (We use the convention of summing over
repeated subscript and superscript indices.) Taking the first and second
derivative of $J^{i_{1}i_{2}}L_{i_{1}}H_{i_{2}}=0$ with respect to $k\,$th
component of $z$ and using the fact that, since $L(z)$ is quadratic,
$L_{i_{1}k_{1}k_{2}}=0$, we have%
\begin{subequations}
\begin{align}
J^{i_{1}i_{2}}L_{i_{1}k_{1}}\mathcal{H}_{i_{2}}+J^{i_{1}i_{2}}L_{i_{1}%
}\mathcal{H}_{i_{2}k_{1}}  &  =0\label{first_derivative_of_bracket}\\
J^{i_{1}i_{2}}L_{i_{1}k_{1}}\mathcal{H}_{i_{2}k_{2}}+J^{i_{1}i_{2}}%
L_{i_{1}k_{2}}\mathcal{H}_{i_{2}k_{1}}+J^{i_{1}i_{2}}L_{i_{1}}\mathcal{H}%
_{i_{2}k_{1}k_{2}}  &  =0. \label{second_derivative_of_bracket}%
\end{align}
Since $\psi=(J\mathcal{H}_{z})^{\top}\mathcal{H}_{zz}(J\mathcal{H}%
_{z})=J^{k_{1}j_{1}}J^{k_{2}j_{2}}\mathcal{H}_{k_{1}k_{2}}\mathcal{H}_{j_{1}%
}\mathcal{H}_{j_{2}}$ we have
\end{subequations}
\begin{align}
\psi_{i_{2}}  &  =J^{k_{1}j_{1}}J^{k_{2}j_{2}}\mathcal{H}_{k_{1}k_{2}i_{2}%
}\mathcal{H}_{j_{1}}\mathcal{H}_{j_{2}}+J^{k_{1}j_{1}}J^{k_{2}j_{2}%
}\mathcal{H}_{k_{1}k_{2}}\mathcal{H}_{j_{1}i_{2}}\mathcal{H}_{j_{2}%
}\nonumber\\
&  +J^{k_{1}j_{1}}J^{k_{2}j_{2}}\mathcal{H}_{k_{1}k_{2}}\mathcal{H}_{j_{1}%
}\mathcal{H}_{j_{2}i_{2}}\label{psi_deriv}\\
&  =J^{k_{1}j_{1}}J^{k_{2}j_{2}}\mathcal{H}_{k_{1}k_{2}i_{2}}\mathcal{H}%
_{j_{1}}\mathcal{H}_{j_{2}}+2J^{k_{1}j_{1}}J^{k_{2}j_{2}}\mathcal{H}%
_{k_{1}k_{2}}\mathcal{H}_{j_{1}i_{2}}\mathcal{H}_{j_{2}}\nonumber
\end{align}
where we combined the last two terms of (\ref{psi_deriv}) by renaming indices
and using the fact that $\mathcal{H}_{k_{2}k_{1}}=\mathcal{H}_{k_{1}k_{2}}.$
Substituting for $\psi_{i_{2}}$ in the Possion bracket $[L,\psi]=L_{z}^{\top
}J\psi_{z}=J^{i_{1}i_{2}}L_{i_{1}}\psi_{i_{2}}$ we have
\begin{align*}
\lbrack L,\psi]  &  =J^{k_{1}j_{1}}J^{k_{2}j_{2}}\left(  J^{i_{1}i_{2}%
}L_{i_{1}}\mathcal{H}_{k_{1}k_{2}i_{2}}\right)  \mathcal{H}_{j_{1}}%
\mathcal{H}_{j_{2}}\\
&  +2J^{k_{1}j_{1}}J^{k_{2}j_{2}}\left(  J^{i_{1}i_{2}}L_{i_{1}}%
\mathcal{H}_{j_{1}i_{2}}\right)  \mathcal{H}_{k_{1}k_{2}}\mathcal{H}_{j_{2}}.
\end{align*}
\newline Using equations (\ref{first_derivative_of_bracket})
and\ (\ref{second_derivative_of_bracket}) and $\mathcal{H}_{j_{1}i_{2}%
}=\mathcal{H}_{i_{2}j_{1}},\ $ $\mathcal{H}_{k_{1}k_{2}i_{2}}=\mathcal{H}%
_{i_{2}k_{1}k_{2}},$ we have%
\begin{align}
\lbrack L,\psi]  &  =-J^{k_{1}j_{1}}J^{k_{2}j_{2}}\left(  J^{i_{1}i_{2}%
}L_{i_{1}k_{1}}H_{i_{2}k_{2}}+J^{i_{1}i_{2}}L_{i_{1}k_{2}}H_{i_{2}k_{1}%
}\right)  H_{j_{1}}H_{j_{2}}\nonumber\\
&  -2J^{k_{1}j_{1}}J^{k_{2}j_{2}}\left(  J^{i_{1}i_{2}}L_{i_{1}j_{1}}H_{i_{2}%
}\right)  H_{k_{1}k_{2}}H_{j_{2}} \label{last_term}%
\end{align}
The last term of (\ref{last_term}) can be expressed as (\ref{new_last_term})
by rearranging terms, using $J^{k_{1}j_{1}}=-J^{j_{1}k_{1}}$ and renaming
indices as shown below.%
\begin{gather}
-2J^{k_{1}j_{1}}J^{k_{2}j_{2}}\left(  J^{i_{1}i_{2}}L_{i_{1}j_{1}}H_{i_{2}%
}\right)  H_{k_{1}k_{2}}H_{j_{2}}=\nonumber\\
-2J^{i_{1}i_{2}}J^{k_{2}j_{2}}\left(  J^{k_{1}j_{1}}L_{j_{1}i_{1}}%
H_{k_{1}k_{2}}\right)  H_{i_{2}}H_{j_{2}}=\nonumber\\
2J^{i_{1}i_{2}}J^{k_{2}j_{2}}\left(  J^{j_{1}k_{1}}L_{j_{1}i_{1}}H_{k_{1}%
k_{2}}\right)  H_{i_{2}}H_{j_{2}}=\nonumber\\
2J^{k_{1}j_{1}}J^{k_{2}j_{2}}\left(  J^{i_{1}i_{2}}L_{i_{1}k_{1}}H_{i_{2}%
k_{2}}\right)  H_{j_{1}}H_{j_{2}} \label{new_last_term}%
\end{gather}
Replacing the last term in (\ref{last_term}) with (\ref{new_last_term}) and
rearranging we have
\begin{equation}
\lbrack L,\psi]=J^{k_{1}j_{1}}J^{k_{2}j_{2}}J^{i_{1}i_{2}}\left(
L_{i_{1}k_{1}}H_{i_{2}k_{2}}-L_{i_{1}k_{2}}H_{i_{2}k_{1}}\right)  H_{j_{1}%
}H_{j_{2}} \label{skew_symmetry}%
\end{equation}
Skew-symmetry with respect to the indicies $k_{1}$ and $k_{2}$ in
(\ref{skew_symmetry}) implies $[L,\psi]=0.$
\end{proof}

\begin{theorem}
[Quadratic Conservation Laws]\label{quad_conservation}Assume $L(z)$ is a
quadratic function and $[L,\mathcal{H}]$ is identically equal to zero. Then
$L(z)$ is exactly conserved at the vertices of a regularized DTH trajectory.
\end{theorem}

\begin{proof}
Since $L(z)$ is quadratic,
$L(z_{k+1})-L(z_{k})=L_{z}(\overline{z}_{k})^{\top }\Delta z_{k}.$
From equation (\ref{reg_DTH_eq_1}) and Lemma \ref{psi_bracket}
we have%
\begin{align*}
L(z_{k+1})-L(z_{k})  &  =L_{z}(\overline{z}_{k})^{\top}\left[  \lambda
_{k}J\mathcal{H}_{z}(\overline{z}_{k})+\mu_{k}J\psi_{z}(\overline{z}%
_{k})\right] \\
&  =\left.  \left(  \lambda_{k}[L,\mathcal{H]+}\mu_{k}[L,\mathcal{\psi
]}\right)  \right\vert _{z=\overline{z}_{k}}\\
&  =0.
\end{align*}

\end{proof}

\begin{corollary}
[Quadratic Conservation Laws]If $L(z)$ is quadratic, $L_{t}\ $and $L_{\wp}%
\ $are both zero, and the Poisson bracket $[L,H]$ is identically
equal to zero, then $L(z)$ is exactly conserved at the vertices of
a regularized DTH trajectory.
\end{corollary}

\begin{proof}
Since $[L,\mathcal{H}]=[L,H]+$ $L_{t}\mathcal{H}_{\wp}-L_{\wp}\mathcal{H}_{t}%
$, then $L_{t}=L_{\wp}=0$ implies $[L,\mathcal{H}]=[L,H]\ $and the corollary
follows from Theorem \ref{quad_conservation}.
\end{proof}

\subsection{Coordinate Invariance}

We briefly consider the coordinate invariance of regularized DTH dynamics. In
the lemma below, we show that $\psi(z)$ is coordinate invariant with respect
to linear, symplectic, coordinate transformations.

\begin{lemma}
\label{coord_inv_psi}Let $z=TZ$ be a linear, symplectic, coordinate
transformation between old coordinates $z$ and new coordinates $Z$. Let
$\mathcal{H}(z)$ be a Hamiltonian function expressed in the old coordinates
and $\mathcal{K}(Z)=\mathcal{H}(TZ)$ be the Hamiltonian function expressed in
the new coordinates. Define $\psi^{\mathcal{H}}(z)=(J\mathcal{H}_{z})^{\top
}\mathcal{H}_{zz}(J\mathcal{H}_{z})$ and $\psi^{\mathcal{K}}(Z)=\left(
J\mathcal{K}_{Z}\right)  ^{\top}\mathcal{K}_{ZZ}\left(  J\mathcal{K}%
_{Z}\right)  $. Then $\psi^{\mathcal{K}}(Z)=\psi^{\mathcal{H}}(TZ).$
\end{lemma}

\begin{proof}
Since $\mathcal{K}(Z)=\mathcal{H}(TZ)$ we have $\mathcal{K}_{Z}=T^{\top
}\mathcal{H}_{z}$ and $K_{ZZ}=T^{\top}\mathcal{H}_{zz}T.$
\begin{align*}
\psi^{K}(Z)  &  =\left(  J\mathcal{K}_{Z}\right)  ^{\top}\mathcal{K}%
_{ZZ}\left(  J\mathcal{K}_{Z}\right) \\
&  =\left(  JT^{\top}\mathcal{H}_{z}\right)  ^{\top}\left(  T^{\top
}\mathcal{H}_{zz}T\right)  \left(  JT^{\top}\mathcal{H}_{z}\right) \\
&  =\mathcal{H}_{z}^{\top}\left(  TJ^{\top}T^{\top}\right)  \mathcal{H}%
_{zz}\left(  TJT^{\top}\right)  \mathcal{H}_{z}%
\end{align*}
Since $T$ is symplectic, $TJT^{\top}=T^{\top}JT=J$ and we have
\begin{align*}
\psi^{K}(Z)  &  =\mathcal{H}_{z}^{\top}J^{\top}\mathcal{H}_{zz}J\mathcal{H}%
_{z}\\
&  =(J\mathcal{H}_{z})^{\top}\mathcal{H}_{zz}(J\mathcal{H}_{z})\\
&  =\psi^{\mathcal{H}}(z)\\
&  =\psi^{\mathcal{H}}(TZ).
\end{align*}

\end{proof}

\begin{theorem}
[Linear Coordinate Invariance]\label{coordinate_inv}The regularized DTH
equations are coordinate invariant under linear, symplectic, coordinate transformations.
\end{theorem}

The proof of Theorem \ref{coordinate_inv} parallels the proof given in
\cite{Shibberu-97} for the DTH equations. The only difference is the use of
Lemma \ref{coord_inv_psi} in the proof of Theorem \ref{coordinate_inv}. (We
point out that the Lagrange multiplier $\lambda_{k}$ and the KKT multiplier
$\mu_{k}$ are both coodinate invariant quantities.)

Ge \cite{Ge-91} demonstrated the coordinate invariance of a variety of
symplectic integrators under linear, symplectic coordinate transformations.
Guibout and Bloch \cite{Guibout-04} demonstrate coordinate invariance using
the larger class of linear, symplectic, discrete-time coordinate
transformations. In fact, it should be possible to demonstrate coordinate
invariance using the even larger class of piecewise-linear, continuous,
symplectic coordinate transformations which are consistent with a special
triangulation of extended phase space. A procedure for demonstrating the
coordinate invariance of DTH dynamics using this larger class of coordinate
transformations was described (in a formal sense) in \cite{Shibberu-92}.
Theorem \ref{coordinate_inv} implies this procedure is also valid for
regularized DTH dynamics. In fact, the regularization described in this
article removes the technical difficulty of identifying principle DTH
trajectories described in \cite{Shibberu-92}.

\section{Numerical Results\label{numerics}}

\subsection{An Algorithm for Regularized SEM Integration}

An algorithm for solving the regularized DTH equations (\ref{reg_DTH_eq_1}%
)--(\ref{reg_DTH_eq_5}) is outlined in Figure \ref{algorithm}. In this
section, we will choose $\psi_{k}=0$ in equation (\ref{reg_DTH_eq_3}). Before
explaining the algorithm in detail, it would be useful to review the simpler
algorithm developed in \cite{Shibberu-92} for solving the DTH equations
(\ref{DTH_equations_1})--(\ref{DTH_equations_2}).

Equations (\ref{DTH_equations_1})--(\ref{DTH_equations_2}) are poorly
conditioned for small time steps $\lambda_{k}$. A direct application of
Newton's method is likely to result in poor convergence. Instead, nested,
Newton iterations are used. The function $\overline{z}_{k}=\overline
{z}(\lambda_{k},z_{k})$ implicity defined by equation (\ref{DTH_equations_1}),
is evaluated using an inner iteration. An outer iteration solves the equation
$g(\lambda)=\mathcal{H}(\overline{z}(\lambda,z_{k}))=0$ for $\lambda_{k}$.
Quadratic convergence of the iterations is proved in \cite{Shibberu-92} for
$\psi(z_{k})\neq0.$ The outer Newton iteration exhibits poor convergence near
$\psi=0.$

The algorithm outlined in Figure \ref{algorithm} also uses nested, Newton
iterations. Near $\psi=0$ however, a bracketed, root-finding procedure is used
in place of the outer Newton iteration. The algorithm is further complicated
by the problem of identifying when $\psi=0$ has been crossed. The procedural
logic needed to compute a DTH trajectory crossing $\psi=0$ appears to be
complex. We now give a more detailed explanation of the algorithm outlined in
Figure \ref{algorithm}.%

%TCIMACRO{\TeXButton{%
%\begin{figure}%
%}{\begin{figure}}}%
%BeginExpansion
\begin{figure}%
%EndExpansion

1.$\qquad\qquad\left\{
\begin{tabular}
[c]{l}%
input $z_{0},\ $set $k=0$\\
solve $\psi(\overline{z}(\lambda,z_{k}))\mathcal{H}(\overline{z}(\lambda
,z_{k}))=0$ for $\lambda\geq0$ to determine $\lambda_{k}$\\
$\overline{z}_{k}=\overline{z}(\lambda_{k},z_{k}),$ $z_{k+1}=2\overline{z}%
_{k}-z_{k}$%
\end{tabular}
\ \right.  $

2.\qquad$%
\begin{tabular}
[c]{l}%
repeat while $k\leq N$%
\end{tabular}
\ $

3.$\qquad\qquad\left\{
\begin{tabular}
[c]{l}%
while $\psi(z_{k})\psi(z_{k+1})>0$ and $k\leq N$\\
\hspace{0.25in}$k=k+1$\\
\hspace{0.25in}solve $\psi(\overline{z}(\lambda,z_{k}))\mathcal{H}%
(\overline{z}(\lambda,z_{k}))=0$ for $\lambda\geq0$ to determine $\lambda_{k}%
$\\
\hspace{0.25in}$\overline{z}_{k}=\overline{z}(\lambda_{k},z_{k}),$
$z_{k+1}=2\overline{z}_{k}-z_{k}$\\
end
\end{tabular}
\ \right.  $

4.\qquad\qquad$\left\{
\begin{tabular}
[c]{l}%
solve $\psi(\overline{z}(\lambda,z_{k}))=0$ to determine $\lambda_{\psi}$\\
$\overline{z}_{\psi}=\overline{z}(\lambda_{\psi},z_{k})$%
\end{tabular}
\right.  $

5.$\qquad\qquad\left\{
\begin{tabular}
[c]{l}%
while $\lambda_{\psi}\geq0\ $and $\mathcal{H}(z_{k})\mathcal{H}(\overline
{z}_{\psi})\leq0$ and $k\leq N$\\
\hspace{0.25in}solve $\mathcal{H}(\overline{z}(\lambda,z_{k}))=0$ for
$0\leq\lambda\leq\lambda_{\psi}$ to determine $\lambda_{k}$\\
\hspace{0.25in}$\overline{z}_{k}=\overline{z}(\lambda_{k},z_{k}),$
$z_{k+1}=2\overline{z}_{k}-z_{k}$\\
\hspace{0.25in}$k=k+1$\\
\hspace{0.25in}solve $\psi(\overline{z}(\lambda,z_{k}))=0$ to determine
$\lambda_{\psi}$\\
\hspace{0.25in}$\overline{z}_{\psi}=\overline{z}(\lambda_{\psi},z_{k})$\\
end
\end{tabular}
\ \right.  $

6.\qquad\qquad$\left\{
\begin{tabular}
[c]{l}%
if ghost trajectory and $\mathcal{H}(z_{k-1})\mathcal{H}(z_{k})>0$\\
\hspace{0.25in}$k=k-1$\\
\hspace{0.25in}solve $\psi(\overline{z}(\lambda,z_{k}))=0$ to determine
$\lambda_{\psi}$\\
\hspace{0.25in}solve $\mathcal{H}(\overline{z}(\lambda,z_{k}))=0$ for
$\lambda\geq\lambda_{\psi}$ to determine $\lambda_{k}$\\
$\hspace{0.25in}\overline{z}_{k}=\overline{z}(\lambda_{k},z_{k}),$
$z_{k+1}=2\overline{z}_{k}-z_{k}$\\
\hspace{0.25in}$k=k+1$\\
end
\end{tabular}
\right.  $

7.\qquad\qquad$\left\{
\begin{tabular}
[c]{l}%
if regularized trajectory\\
\hspace{0.25in}solve $\left\{
\begin{tabular}
[c]{r}%
$\mathcal{H}(\overline{z}(\lambda,\mu,z_{k}))=0$\\
$\psi(\overline{z}(\lambda,\mu,z_{k}))=0$%
\end{tabular}
\ \ \right.  $ to determine $\lambda_{k}$ and $\mu_{k}$\\
\hspace{0.25in}$\overline{z}_{k}=\overline{z}(\lambda_{k},\mu_{k},z_{k}),$
$z_{k+1}=2\overline{z}_{k}-z_{k}$\\
\hspace{0.25in}$k=k+1$\\
\hspace{0.25in}solve $\psi(\overline{z}(\lambda,z_{k}))=0$ to determine
$\lambda_{\psi}$\\
\hspace{0.25in}solve $\mathcal{H}(\overline{z}(\lambda,z_{k}))=0$ for
$\lambda\geq\max(0,\lambda_{\psi})$ to determine $\lambda_{k}$\\
\hspace{0.25in}$\overline{z}_{k}=\overline{z}(\lambda_{k},z_{k}),$
$z_{k+1}=2\overline{z}_{k}-z_{k}$\\
\hspace{0.25in}$k=k+1$\\
end
\end{tabular}
\right.  $

8.\qquad\qquad$\left\{
\begin{tabular}
[c]{l}%
solve $\mathcal{H}(\overline{z}(\lambda,z_{k}))=0$ for $\lambda\geq0$ to
determine $\lambda_{k}$\\
$\overline{z}_{k}=\overline{z}(\lambda_{k},z_{k}),\ z_{k+1}=2\overline{z}%
_{k}-z_{k}$\\
$k=k+1$%
\end{tabular}
\right.  $

\ \ \qquad$%
\begin{tabular}
[c]{l}%
end
\end{tabular}
$%

%TCIMACRO{\TeXButton{\caption}{\caption
%{Algorithm for computing ghost and regularized DTH trajectories.}}}%
%BeginExpansion
\caption{Algorithm for computing ghost and regularized DTH trajectories.}%
%EndExpansion%
%TCIMACRO{\TeXButton{\label}{\label{algorithm}}}%
%BeginExpansion
\label{algorithm}%
%EndExpansion
%

%TCIMACRO{\TeXButton{%
%\end{figure}%
%}{\end{figure}}}%
%BeginExpansion
\end{figure}%
%EndExpansion

The function $\overline{z}(\lambda_{k},\mu_{k},z_{k})$ implicitly defined by
equation (\ref{reg_DTH_eq_1}) is evaluated using Newton's method. When
$\psi(\overline{z}_{k})$ $\neq0,$ equation (\ref{reg_DTH_eq_4}) implies
$\mu_{k}=0.$ We use the abbreviation $\overline{z}(\lambda_{k},z_{k})$ for
$\overline{z}(\lambda_{k},\mu_{k},z_{k})$ when $\mu_{k}=0$.

The first block in Figure \ref{algorithm} initializes the algorithm. (We
assume $\psi(\overline{z}_{0})\neq0.$) The value of $\wp_{0}$ (the momentum
conjugate to time) determines the initial time step $\lambda_{0}.$ If $\wp
_{0}\ $is chosen so that $\mathcal{H}(z_{0})=0,$ then $\lambda_{0}=0.$
Therefore, a value for $\wp_{0}$ should be chosen so that $\mathcal{H}(z_{0})$
is sufficiently small but nonzero.

Vertex points $z_{k}$, $k=2,\ldots,N,$ are computed within block 2. To avoid
ill-conditioning of the equation $\mathcal{H}(\overline{z}(\lambda,z_{k}))=0$
near $\psi=0,$ we solve the equation $\psi(\overline{z}(\lambda,z_{k}%
))\mathcal{H}(\overline{z}(\lambda,z_{k}))=0$. Linear segments which do not
cross $\psi=0$ are computed in block 3. If $\psi(\overline{z}(\lambda
_{k},z_{k}))=0,$ the condition $\psi(z_{k})\psi(z_{k+1})\leq0$ indicates
$\psi=0$ has been crossed and the algorithm enters block 4 where
$\lambda_{\psi}$ and $\overline{z}_{\psi}$ are computed and used in block 5.
In block 5, a bracketed root-finding procedure is used to solve the now
ill-conditioned equation $\mathcal{H}(\overline{z}(\lambda,z_{k}))=0$. If a
bracket can not be found, the algorithm enters either block 6 and computes a
ghost trajectory or block 7 and computes a regularized trajectory. Block 8 is
need to prevent the computation of trajectories which immediately recross
$\psi=0.$

\begin{figure}
\begin{center}
\subfigure[Parameter]{\label{one_parm}\includegraphics
[height=1in]{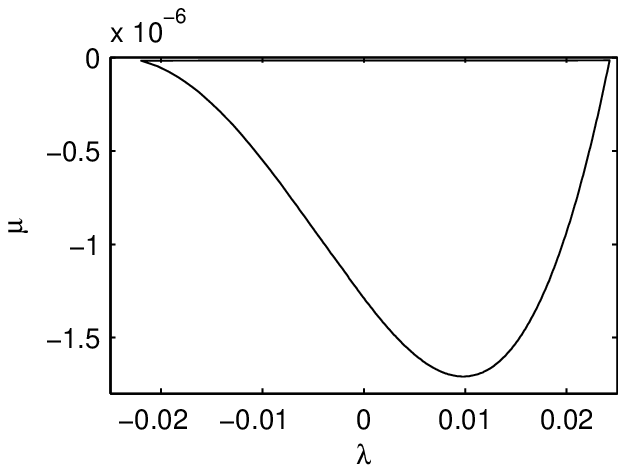}}
\subfigure[$\lambda<$ 0,
$\mu$-small]{\includegraphics[height=1in]{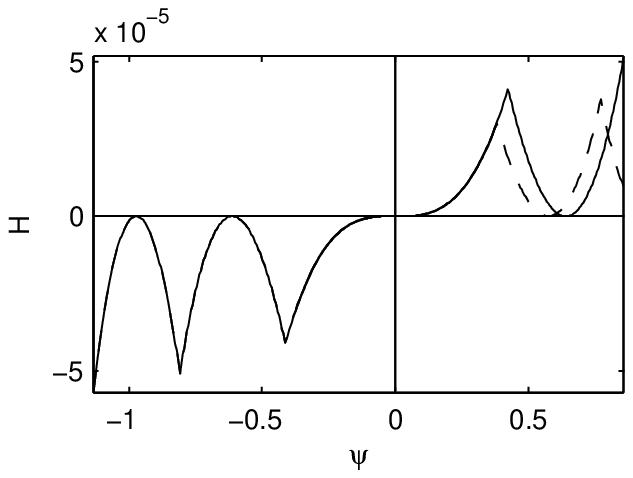}}
\subfigure[$\lambda<$
0]{\includegraphics[height=1in]{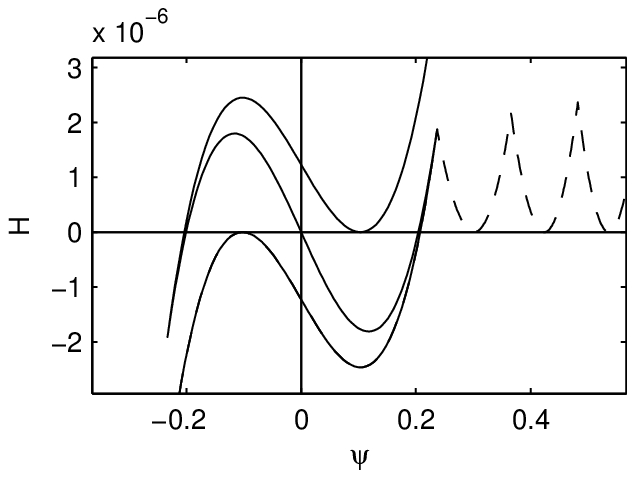}}
\subfigure[$\lambda=$
0]{\includegraphics[height=1in]{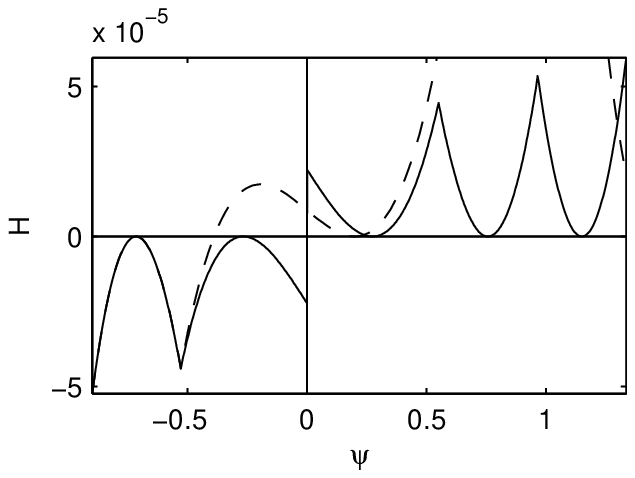}}
\subfigure[$\lambda>$
0]{\includegraphics[height=1in]{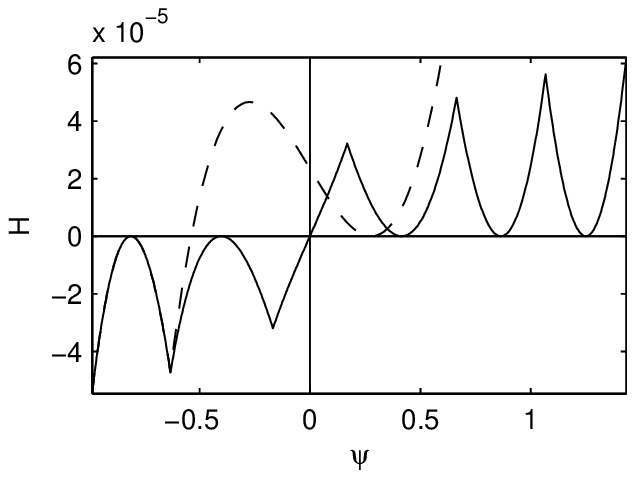}}
\subfigure[$\lambda>$ 0,
$\mu$-small]{\includegraphics[height=1in]{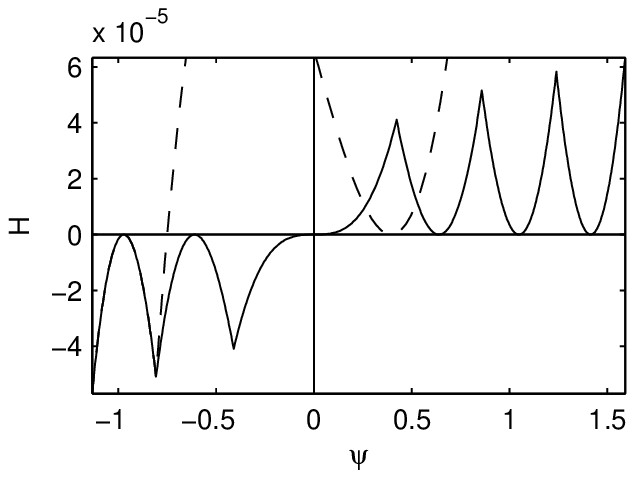}}
\end{center}
\caption{A one parameter family of DTH trajectories crossing $\psi
= 0$. Dashed curves are ghost trajectories. Solid curves are
regularized trajectories.} \label{crossings}
\end{figure}

Blocks 7 and 8 take into account the different ways DTH trajectories can cross
$\psi=0.$ These blocks are best understood after viewing an animation of a
one-parameter family of $\psi=0$ crossings. Snapshots of this animation are
given in Figure \ref{crossings}.%

%TCIMACRO{\TeXButton{Figure: crossings}{\begin{figure}
%\begin{center}
%\subfigure[Parameter]{\label{one_parm}\includegraphics
%[height=1in]{one_parm_family.eps}}
%\subfigure[$\lambda<$ 0,
%$\mu$-small]{\includegraphics[height=1in]{psi_crossing_1.eps}}
%\subfigure[$\lambda<$
%0]{\includegraphics[height=1in]{psi_crossing_2.eps}}
%\subfigure[$\lambda=$
%0]{\includegraphics[height=1in]{psi_crossing_3.eps}}
%\subfigure[$\lambda>$
%0]{\includegraphics[height=1in]{psi_crossing_4.eps}}
%\subfigure[$\lambda>$ 0,
%$\mu$-small]{\includegraphics[height=1in]{psi_crossing_5.eps}}
%\end{center}
%\caption{A one parameter family of DTH trajectories crossing $\psi
%= 0$. Dashed curves are ghost trajectories. Solid curves are
%regularized trajectories.} \label{crossings}
%\end{figure}}}%
%BeginExpansion
%EndExpansion

Murua \cite{Murua-97} has developed an efficient iteration which avoids the
nested iterations use to solve $\mathcal{H}(\overline{z}(\lambda,z_{k}))=0$.
(Murua has also developed an iteration which does not require evaluation of
the Hessian matrix of the Hamiltonian function.) It is likely that the
algorithm outlined in Figure \ref{algorithm} could be made significantly more
efficient by using Murua's iteration. The author has take a first step in this
direction by modifying Murua's iteration so that it can be used to compute
regularized segments crossing $\psi=0.$

\subsection{Qualitative Behavior of Regularized SEM Integation}

Numerical computations for the pendulum and Kepler's one body
problem confirm that the regularization described in this article
conserves SEM properties. The energy (Hamiltonian) is conserved to
roundoff error at midpoints of DTH trajectories. Angular momentum
is conserved to roundoff error at vertices for Kepler's problem in
Cartesian coordinates. Symplecticness is
verified by computing the derivative $dz_{N}/dz_{0}$ of the map $z_{N}%
(z_{0}).$ The matrix $\left(  dz_{N}/dz_{0}\right)  ^{\top}J\left(
dz_{N}/dz_{0}\right)  \ $is found to equal $J$ to roundoff error.
Time-reversibility is also confirmed to hold to roundoff error. Numerical
computations quantifying the accuracy and efficiency of the regularization
have not yet been completed.

One of the peculiarities observed in regularized SEM integration is the
occurrence of negative time steps. Negative time steps violate the
monotonic-increasing property of time. The DTH trajectories become
multi-valued functions of time. Lee \cite{Lee-87} foresaw this possibility and
suggested the remedy of relinking vertices to maintain the monotonicity of time.

\begin{figure}
\begin{center}
\subfigure[$q=0,p=3$]{\label{}\includegraphics
[height=1.5in]{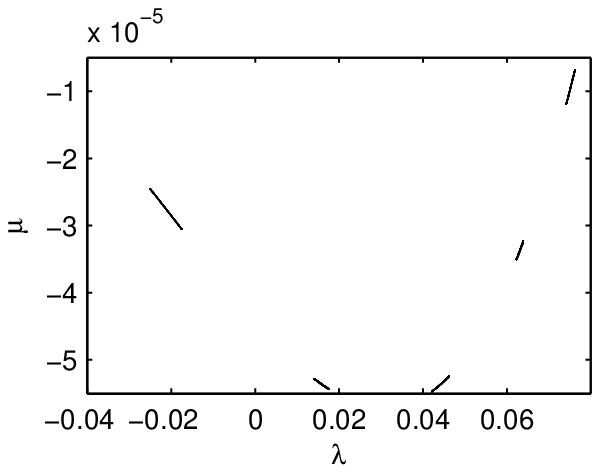}}
\subfigure[$q=0,p=2.2$]{\label{}\includegraphics[height=1.5in]{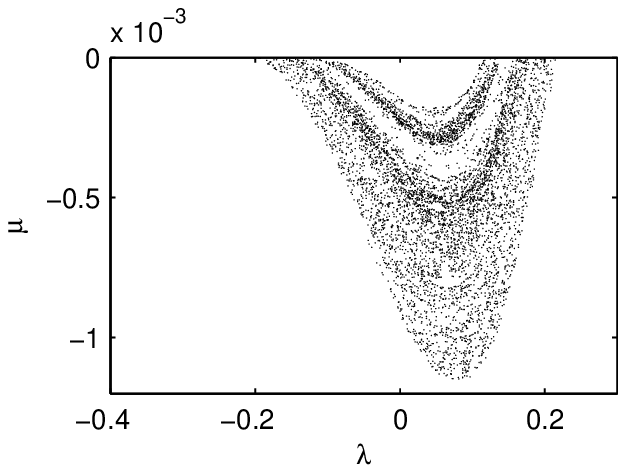}}
\end{center}
\caption{Behavior of $\lambda$ vs $\mu$ when $\psi=0$ for a single
regularized DTH trjectory of the nonlinear pendulum.} \label{chaos}
\end{figure}

Finally, another peculiarity observed for regularized SEM
integration, for the case $\psi_{k}=0$ in equation
(\ref{reg_DTH_eq_3}), is apparently chaotic behavior near the
separatrix of the pendulum. (See Figure \ref{chaos}.) It may be
possible to regularize DTH dynamics further by choosing nonzero
values for $\psi_{k}.$%

%TCIMACRO{\TeXButton{Figure: chaos}{\begin{figure}
%\begin{center}
%\subfigure[$q=0,p=2.2$]{\label{}\includegraphics
%[height=1.5in]{lam_vs_mu_1.eps}}
%\subfigure[$q=0,p=3$]{\label{}\includegraphics[height=1.5in]{lam_vs_mu_2.eps}}
%\end{center}
%\caption{Behavior of $\lambda$ vs $\mu$ when $\psi=0$ for a single
%regularized DTH trjectory of the nonlinear pendulum.} \label{chaos}
%\end{figure}
%}}%
%BeginExpansion
%EndExpansion

\bibliographystyle{hplain}
\bibliography{DTH_dynamics}

\end{document}